\documentclass{amsart}
\usepackage{lineno,hyperref}
\usepackage{amsfonts,amssymb,amsmath}
\usepackage{graphicx}
\usepackage[usenames,dvipsnames]{color}
\usepackage[all]{xy}
\usepackage{tikz}
\usepackage[normalem]{ulem}

\usepackage{amsthm}
\usepackage{amssymb}
\usepackage{amscd}
\usepackage{amsfonts}
\usepackage{amsbsy}
\usepackage{enumerate}
\usepackage{graphicx}
\usepackage[dvips]{psfrag}
\usepackage{subfigure}

\newtheorem{theorem}{Theorem}

\newtheorem{proposition}{Proposition}

\newtheorem{remark}{Remark}

\graphicspath{{./figs/}}
\DeclareGraphicsExtensions{.pdf,.eps,.png,.jpg}

\begin{document}

\title[Hopf and Bautin bifurcations in a 3D model]{Hopf and Bautin bifurcations in a 3D model for pest leafhopper with stage structure and generalist predatory mite} 


\author[Alvarez-Ram\'{\i}rez]{Martha Alvarez-Ram\'{\i}rez }
\address{Dept. de Matem\'aticas, UAM--Iztapalapa,
09340 Iztapalapa,  Ciudad de M\'exico,  M\'exico}
\email{mar@xanum.uam.mx}

\author[Garc\'{\i}a Rivera]{Marco Polo Garc\'{\i}a Rivera}
\email{riveramarcopolo@hotmail.com}

\author[Ortiz Santosz]{Ahida Ortiz Santos}
 \email{ahidaosz@gmail.com}

\begin{abstract}
In a recent paper of Yuan and  Zhu (J. Differential Equations 321(2022) 99-129), the nature of dynamics of  
a generalist predator and prey with stage structure  is modeled as a three-dimensional coupled nonlinear differential system. A detailed analysis of the bifurcation shows that the model exhibits a high complexity in its dynamics, which arises from the use of predatory mites as agent for controlling stage structures of tea green leafhopper pest. Unfortunately,  there is a mistake in Proposition 2.2, item (3), also  the hypothesis in Proposition 3.3  is incorrect.  
In this paper, we revisit the model and straighten those mentioned errors as  in item (3),  changed the existing hypothesis  by a suitable one and give a corrected proof of  Proposition 3.3 in its new form.  
Also, we show   that the model does undergo both Hopf bifurcation and Bautin bifurcation, and numerical examples are given to support the analytic results.
\end{abstract}



\maketitle




\section{Introduction}
The production of tea is frequently threatened by one of the most prevalent insect plagues insect plagues  in many Asian nations is the tea green leafhopper or {\em Empoasca onukii}, and the mite {\em Anystis baccarum} has been  applied as a plague control agent, which is a predatory that is beneficial to the tea plantations. For the purpose of researching the interactions between adult populations, Yuan and Zhu \cite{Yuan2022} constructed a three-dimensional ordinary differential equation model, which  enables them to use bifurcation theory to investigate the causes and reasons existing behind their complex dynamics.
However, in contrast to their earlier paper \cite{Yuan2021}, this takes into account the interaction between generalist predator mite and two different developmental stages of {\em E. onukii}.

In  \cite{Yuan2022}, the authors consider a  model to analyze dynamics and inform biological control.
In particular, they  perform a bifurcation study of a 3D nonlinear system, and found some complex bifurcation phenomena, such as Bogdanov-Takens bifurcation, Hopf bifurcation,  saddle-node bifurcation of codimension 1 and 2, and bifurcations of nilpotent singularities of elliptic and focus type of codimension 3.  In addition, they found that  all  codimension 3 bifurcations in the two-dimensional center manifold of the system are connected by the nilpotent focus of codimension 4, which acts as an organizing center.
However, there is a mistake in Proposition 3.3 where the authors assume that $\kappa_2 < \widehat{\kappa}_2$ must be true in order to show that a Hopf bifurcation  exists at the coexistence equilibrium points, this hypothesis is incorrect. 

Motivated by the aforementioned research, our work is devoted to  fix those mentioned gaps and some others errors.
According to the known results,  the existence of a  Bautin bifurcation (generalized Hopf) in the model described by Yuan and Zhu \cite{Yuan2022} was not supported by any evidence.
Thus, in this manuscript we will present numerical data to support the existence of Hopf bifurcation  of codimension 1 and  2 (Bautin).

In order to obtain results that, up to this moment,  cannot be obtained analytically, as well as to generate phase portraits and bifurcation plots, we heavily rely on advanced continuation and bifurcation techniques implemented in the {\em MatCont} software package \cite{matcont}.
Finally, we stress that some calculations  were made easier by using  the computer algebra system Wolfram  {\em Mathematica}, version 12.

We want to emphasize how some calculations were made easier by using the computer algebra system.

The paper is organized as follows. In section \ref{sec1} we introduce the model system.
 In order to make this work self-contained and facilitate the lecture to the reader, we summarize briefly some current known facts 
regarding  equilibrium points  in Section \ref{sec_pts}.
The local analysis of the  stability of the interior (coexistence) equilibrium points is provided in Section \ref{sec_sta}.
In Section \ref{sec_hopf}  Hopf bifurcation of the interior equilibrium point of the system is discussed. 
To illustrate the theoretical findings in, we used phase portraits and numerical bifurcation diagrams
 that are given in Section \ref{sec_nume}.  Furthermore, we use the numerical software MatCont \cite{matcont} to compute curves of equilibria and to compute several bifurcation curves. We particularly approximate a family of limit cycles emanating from a Hopf point. Our obtained results give a step forward  improving some results given in \cite{Yuan2022}.

\section{Mathematical model formulation}\label{sec1}
In the following, we briefly discuss the underlying model considering the tea plantation that has a designated area for growing tea trees.To introduce the mathematical model, we use $E_1$ and $E_2$ to represent, respectively, the  {\em E. onukii} 
hatchlings (adults and nymphs) as well as eggs  at time $t$, while $M(t)$ denotes the  {\em A. baccarum}.  
By employing a modeling approach that was used in   \cite{Yuan2021}, the model with 
generalist predator-prey and  functional response, Holling type II,   also with a prey life cycle stage structure  of prey is
 described by the differential equations system
 \begin{equation}\label{eq_model}
 \begin{array}{l}
 \dot{E}_1 = r_1E_2-\alpha E_1, \vspace{0.2cm}\\
 \dot{E}_2 = \alpha E_1- \kappa_1 E_2^2 - \dfrac{m E_2 M}{a+E_2}, \vspace{0.2cm}\\
  \dot{M} = r_2 M -\kappa_2 M^2 + \dfrac{cm E_2 M}{a+E_2},
 \end{array}
 \end{equation}
where the dot means derivative with respect to the time $t$.  All the parameters are positive and their  biological meanings are the following: 
The rate at which{\em E. onukii} eggs hatch is  $\alpha$, $a$ is the number of {\em E. onukii} that {\em A. baccarum} has managed to locate and successfully capture in an average amount of time,
the conversion rate is given by $c$, the maximum number of {\em E. onukii} that {\em A. baccarum} can sustain during the foraging period is given by $m$, the adult {\em E. onukii} rate of oviposition is given by $r_1$, and the constant intrinsic growth rate of {\em A. baccarum} is indicated by $r_2$, $\kappa_1$ and $\kappa_2$ are  the intra-specific levels of  competition among adult of 
{\em E. onukii} and {\em A. baccarum}, respectively. 

It is important to remark that $a>0$ indicates that the predator is generalist and that if it does not exist available prey 
it has an  alternative  source of food.

Before beginning our search, we shall reduce the number of parameters  in model  \eqref{eq_model}. To do so,
firstly, we introduce a change of variables and a time rescaling  $(E_1,E_2,M,t)\to (aX_1,aX_2,acY,\tau)$, given by
$$
E_1 = a X_1, \quad \qquad  E_2 = a X_2, \quad \qquad M = a c Y, \quad\qquad t= \dfrac{1}{cm}\tau.
$$
Observe that we still denote $\tau$ by $t$ and  we have renamed the variables $(X_1,X_2,Y)$ again as $(E_1,E_2,M)$
for convenience.  Thus,  system \eqref{eq_model} can be reduced as follows 
 \begin{equation}\label{eq_model2}
 \begin{array}{l}
 \dot{E}_1 = \hat{r}_1E_2 - \hat{\alpha}E_1, \vspace{0.2cm}\\
 \dot{E}_2 = \hat{\alpha} E_1- \hat{\kappa}_1E_2^2 - \dfrac{E_2 M}{1+E_2}, \vspace{0.2cm}\\
  \dot{M} = \hat{r}_1 M -\hat{\kappa}_2 M^2 + \dfrac{E_2 M}{1+E2},
 \end{array}
 \end{equation}
where $\hat{r}_1=\frac{r_1}{cm}$, $\hat{\alpha}= \frac{\alpha}{cm}$, $\hat{\kappa}_1=\frac{a\kappa_1}{cm}$, $\hat{\kappa}_2=\frac{a\kappa_2}{cm}$ and dropping all the hats over the parameters, the system becomes
  \begin{equation}\label{eq_model3}
 \begin{array}{l}
 \dot{E}_1 = r_1E_2 - \alpha E_1, \vspace{0.2cm}\\
 \dot{E}_2 = \alpha E_1- \kappa_1 E_2^2 - \dfrac{E_2 M}{1+E2}, \vspace{0.2cm}\\
  \dot{M} = r_1 M - \kappa_2  M^2 + \dfrac{E_2 M}{1+E_2}.
 \end{array}
 \end{equation}
 Hence, this system has only five free parameters, namely   $\alpha$, $r_1$, $r_2$, $\kappa_1$, $\kappa_2$. These will be used as  bifurcation parameters. 
 The reader should be warned that system \eqref{eq_model3} is topologically equivalent to system \eqref{eq_model} except at the singularity $E_2=-a$.

\section{Existence of equilibria}\label{sec_pts}
In this section we will summarize the results on feasible equilibrium points of the model system  \eqref{eq_model3}
obtained in \cite{Yuan2022}, where 
it was shown that model \eqref{eq_model3}  displays up to three coexistence equilibria and three boundary equilibria in the first octant.

 In the paragraphs below we  summarize the main findings obtained  there.
 
The origin $S_0=(0,0,0)$ is always an equilibrium point. Additionally,  $S_{10}=(p_1K_1,K_1,0)$ with $K_1=\dfrac{r_1}{\kappa_1}$ and $p_1=\dfrac{r_1}{\alpha}$ is a non-negative predator-extinction equilibrium, which refers to the proportion of eggs to newly hatched (nymphs and adults) individuals of  the {\em E. onukiiat} equilibrium. Similar to this, the pest-free equilibrium $S_{01}= (0,0,K_2)$ exists with $K_2=\dfrac{r_2}{\kappa_2}$. 
Here,  $K_1$ is the maximum capacity for {\em E. onukii} (nymphs and adults) without of {\em A. baccarum}, while $K_2$ is the  loading capacity for {\em A. baccarum} when  {\em E. onukii} is not present, for the specific tea plantation.
Concerning the possible strictly positive equilibrium points corresponding to points  located at the intersection of the nullclines
$\dot{E}_1=0$, $\dot{E}_2=0$ and $\dot{M}=0$ in the first octant of the phase space we have that if these exist, they indicate
the coexistence of the three species in equilibrium. This is denoted by  $\overline{S}=(\overline{E}_1, \overline{E}_2,\overline{M})$.
Solving the system of nullclines equations  for \eqref{eq_model3}, we get  that 
the  intersection points are given by the following  curves in the $E_2M$ plane: 
\begin{align}
& U_1(E_2)= M= \frac{(\alpha E_1-\kappa_1 E_2^2)(a+E_2)}{mE_2},\\
& U_2(E_2)= M= \frac{1}{\kappa_2} \left(r_2 +\dfrac{c m E_2}{a+E2}\right),
\end{align}
where $E_1=p_1E_2$.  Equating $U_1(E_2)$ and $U_2(E_2)$, 
after some computations  we get  that $E_2$  is determined by the solution of the equation
\begin{equation}\label{U3}
U(E_2) = E_2^3 -  \left(\dfrac{r_1}{\kappa_1} - 2 a\right)E_2^2+ \left( a^2-2a \dfrac{r_1}{\kappa_1}
 + \dfrac{m(cm+r_2)}{\kappa_1\kappa_2}\right)E_2
+\dfrac{a^2r_1}{\kappa_1\kappa_2}\left( \dfrac{mr_2}{ar_1}-\kappa_2\right)=0,
\end{equation}
which may have up to three positive roots. 
Therefore, the model system \eqref{eq_model3} can have maximum three feasible interior equilibrium
points. By some trivial computations, we  get  that $E_2=0$  associated with the pest-free equilibrium $S_{01}=(0,0,K_2)$ is a solution of 
\eqref{U3}, which requires 
 \begin{equation}\label{kappa2b}
 \kappa_2 = \overline{\kappa}_2 = \dfrac{r_2m}{r_1a}.
\end{equation}
 
 Due to the difficulty to determine the exact solutions of equation \eqref{U3}, we look closely at the Cardano's formula  for the third
degree algebra equation  provides a criterion,  whether  one, two or three real solution of \eqref{U3} exist. This involves a {\em discriminant} whose sign determines 
 the nature of roots. This can be expressed as follows:
 $$
 \Delta_0= \Delta_0(\kappa_2)= \dfrac{m^2}{\kappa_1^4\kappa_2^3}\big(a_2(\kappa_1)\kappa_2^2+ a_1(\kappa_1)\kappa_2+ a_0(\kappa_1) \big),
 $$
 where 
 \begin{align*}
 & a_2(\kappa_1)= 4ac(a\kappa_1+r_1)^3,\\
 & a_1(\kappa_1)=a^2(8c^2m^2-20cmr_2-r_2^2)\kappa_1^2- 2ar_1\delta_1(10cm+r_2)\kappa_1-r_1^2\delta_1^2,\\
 & a_0(\kappa_1)= 4m\delta_1^3\kappa_1,
 \end{align*}
 with $\delta_1= r_2+cm$. In addition, the discriminant with respect to $\kappa_2$
associated to the second order polynomial $f(\kappa_1)= a_2(\kappa_1)\kappa_2^2+ a_1(\kappa_1)\kappa_2+ a_0(\kappa_1)$
 becomes
 $$ \Delta_1 (\kappa_1)= (a\kappa_1r_2+r_1\delta_1)(r_1\delta_1-a\kappa_1\delta_2)^3, $$
 with $\delta_2= 8 cm -r_2$.  Let 
 $$ \kappa_1^*= \dfrac{r_1\delta_1}{a\delta_2}$$
 the unique positive root of  $\Delta_1 (\kappa_1)=0$. Therefore, when $0< \kappa_1 < \kappa_1^*$, $\Delta_0 (\kappa_2)=0$ can have two positive roots, 
 denoted by $\kappa_2^-$ and $\kappa_2^+$. 
 Furthermore,  if $\kappa_1^*$ into $\Delta_0 (\kappa_2)=0$, then one gets
 $$ \kappa_2^*= \dfrac{\delta_1^2\delta_2}{27r_1ac^2m},$$ and follows that 
 $\Delta_0 (\overline{\kappa}_2)=0$ at 
 $$
 \overline{\kappa}_1= \dfrac{r_1(r_2-cm)}{r_2a}, \qquad \text{with} \quad \textcolor{red}{r_2} > cm.
 $$
 Therefore, we conclude the following result.
 
\begin{figure}[hpbt]
\centering
\includegraphics[scale=0.4]{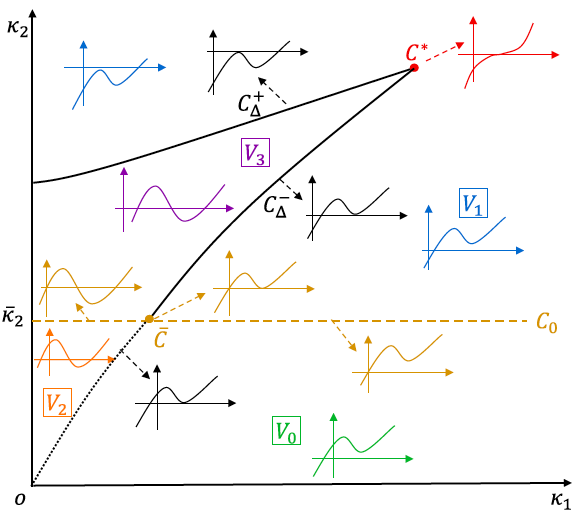}
\caption{Existence of positive equilibriums (positions and number) for the system \eqref{eq_model3} with parameters $\kappa_1$ and $\kappa_2$.
The figure was taken from \cite{Yuan2022}.} \label{fig_region}
\end{figure}

\begin{proposition}[\cite{Yuan2022}]
In the $(\kappa_1, \kappa_2)$ plane, three curves
\begin{align*}
C_0: \kappa_2 &= \overline{\kappa}_2, \qquad \kappa_1 > 0,\\
C_{\triangle}^- :\kappa_2 &= \kappa_2^-(\kappa), \qquad 0< \kappa_1 < \kappa^*, \quad \kappa_2>0, \\
C_{\triangle}^+ :\kappa_2 &= \kappa_2^+(\kappa), \qquad 0< \kappa_1 < \kappa^*, \quad  \kappa_2>0,
\end{align*}
divide the region $\kappa_1>0$, $\kappa_2>0$ into four subregions $V_0$, $V_1$, $V_2$ y $V_3$ {\rm (}see Figure  \ref{fig_region}{\rm )}:
\begin{align*}
V_0 &= \{ (\kappa_1, \kappa_2) \mid  \kappa_1 >0, \quad (\kappa_2 < \kappa^-(\kappa_1)) \cap (\kappa_2 < \overline{\kappa}_2) \}, \\
V_1 &= \{(\kappa_1, \kappa_2)\mid \kappa_1 >0, \quad  \kappa_2 < \bar{\kappa_2} \} \setminus V_3, \\
V_2 &= \{(\kappa_1, \kappa_2)\mid \kappa_1 >0, \quad \kappa_2^-(\kappa_1) < \kappa_2 < \overline{\kappa}_2 \} , \\
V_3 &= \{(\kappa_1, \kappa_2) \mid 0 < \kappa_1 < \kappa_1^*,  \quad \kappa_2^-(\kappa_1) < \kappa_2 < \kappa_2^+ (\kappa_1) \cap (\kappa_2 > \overline{\kappa}_2) \}.
\end{align*}
\begin{enumerate}
\item along $C_0$, where  $\overline{C} = (\overline{\kappa}_1, \overline{\kappa}_2)$
\begin{itemize}
\item if $\kappa_1 > \overline{\kappa}_1$, there is no coexistence equilibria;
\item if $\kappa_1 = \overline{\kappa}_1$, there is a coexistence equilibrium of multiplicity 2, $S_{23}$; 
\item if $\kappa_1 < \overline{\kappa}_1$, there are two coexistence equilibria, $S_2$ y $S_3$.
\end{itemize}
\item along  $C_{\triangle}^+$, there is a  coexistence equilibrium $S_3$  of multiplicity 1 and a coexistence equilibrium of multiplicity 2,  $S_{12}$.
\item along $C_{\triangle}^-$, 
\begin{itemize}
\item if $\kappa_2 > \overline{\kappa}_2$, there is a coexistence equilibrium   $S_1$ of multiplicity 1 and a coexistence equilibrium of multiplicity 2, $S_{23}$;
\item if $\kappa_2 \leq \overline{\kappa}_2$,  there is a coexistence equilibrium of multiplicity 2, $S_{23}$;
\end{itemize}
\item at  $C^* (\kappa_1^* , \kappa_2^*) = \left( \frac{r_1 \delta_1}{a \delta_2}, \frac{\delta_1^2 \delta_2}{27 r_1 a c^2 m} \right)$ there is a coexistence equilibrium point of multiplicity 3, $S_{123}= \left( \frac{r_1 a \delta_3}{\delta_1 \alpha}, \frac{a \delta_3}{\delta_1}, \frac{18 r_1 a c^2 m}{\delta_1 \delta_2}  \right)$   where $0< r_2 < 2cm$.
\end{enumerate}

This allows to split the first quadrant of the positive parameter plane $\kappa_1$-$\kappa_2$ into four basic regions, labeled
$V_0$, $V_1$, $V_2$ and $V_3$, where there are 0, 1, 2, and 3 positive equilibrium points, respectively. They are shown in Figure \ref{fig_region}.  Let $S_1$, $S_2$ and $S_3$ be the three  coexistence equilibria,  as long as there are exist, where its coordinate $\overline{E}_2$  satisfies $\overline{E}_2 (S_1) < \overline{E}_2 (S_2) < \overline{E}_2 (S_3)$.
The positive equilibrium where $S_1$ and $S_2$ coalesce is denoted by $S_{12}$, if it exists.
$S_{23}$ and $S_{123}$ may adopt equivalent definitions.

According to the Figure \ref{fig_region}, the system \eqref{eq_model3} 
there is a chance that the system \eqref{eq_model3} has a single 
coexistence equilibrium  $S_1$ or $S_3$, one degenerate coexistence equilibrium $S_{23}$ or $S_{123}$, two coexistence equilibria simultaneously 
$S_1$, $S_{23}$ or $S_{12}$, $S_3$ or $S_2$, $S_3$, or three coexistence equilibria simultaneously $S_1$, $S_2$, and $S_3$.
\end{proposition}

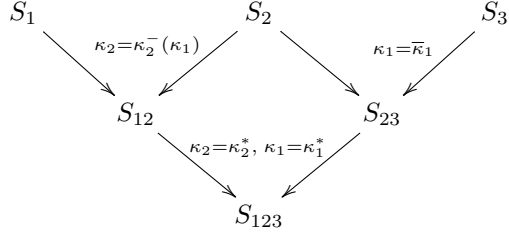
\begin{figure}
\centering
$$ \xymatrix{
S_1\ar[rd]&& S_2\ar[rd] \ar[ld]_{\kappa_2=\kappa_2^-(\kappa_1)}&&S_3\ar[ld]_{\kappa_1= \overline{\kappa}_1}\\
& S_{12} \ar[rd]&& S_{23}\ar[ld] _{\kappa_2 = \kappa_2^*,\;  \kappa_1 = \kappa_1^*} &\\
&& S_{123}    && \\}$$
\caption{The diagram shows the bifurcation of coexistence equilibria. The points $S_1$ and $S_2$ collide when
the parameter  $\kappa_2=\kappa_2^-(\kappa_1)$ and the coinciding equilibria is $S_{12}$.
In a similar way, the equilibria $S_2$ and $S_3$ collide when $\kappa_1= \overline{\kappa}_1$ giving rise to
$S_{23}$.
If $\kappa_1=\kappa_1^*$ and $\kappa_2=\kappa_2^*$  it occurs  that 
all equilibria collapses, that is, 
$S_{12}=S_{23}=S_{123}$.}\label{pts}
\end{figure}

A schematic diagram with a combined collection of  coexistence equilibria and how they are transformed into new ones after some pair of existing or new aquilibria collapse
 is given in Figure \ref{pts}. 
 \section{Stability analysis  of the coexistence equilibria}\label{sec_sta}
 In this section, we shall  discuss the local stability of the model system \eqref{eq_model3} at the coexistence
equilibrium point $\overline{S}$.
The Jacobian matrix of system \eqref{eq_model3}  evaluated 
at  $\overline{S}$ has the form
\begin{equation}\label{matrizJ}
J(\overline{S}) = \left( \begin{array}{ccc} 
-\alpha & r_1 & 0 \\
\alpha & P(\overline{E}_2) U'_1(\overline{E}_2) & -P(\overline{E}_2)\\
0 & \kappa_2 \overline{M} U'_2(\overline{E}_2)  & -\kappa_2  \overline{M}
\end{array}
\right),
\end{equation} 
where prime ($'$) denotes the  derivative with respect to $\overline{E}_2$ and  $P(\overline{E}_2)= \dfrac{m \overline{E}_2}{a+\overline{E}_2}$.

In order to perform an analysis comparative to the one presented in \cite{Yuan2022}, we write 
the {\em right} associated characteristic polynomial of the Jacobian matrix $J(\overline{S})$ as
\begin {equation}\label{p_mex}
p(\lambda)=  -\lambda^3 - A_2 \lambda^2 -A_1 \lambda -A_0
\end{equation}
where
\begin{equation}\label{As}
\begin{array}{l}
A_2 = \alpha +\kappa_2 \overline{M}- P(\overline{E}_2) U'_1(\overline{E}_2),\vspace{0.2cm}\\
A_1 = \alpha \big [\kappa_2 \overline{M} - P(\overline{E}_2) U'_1(\overline{E}_2)  - \textcolor{red}{\bf r_1}\big ] + 
P(\overline{E}_2) \kappa_2 \overline{M}  \big[ U'_2(\overline{E}_2) -U'_1(\overline{E}_2)\big],\vspace{0.2cm}\\
A_0 = \alpha \kappa_2 \overline{M} \big[ P(\overline{E}_2) \big( U'_1(\overline{E}_2) -U'_2(\overline{E}_2)\big) +r_1\big].
\end{array}
\end{equation}

By now, the reader should be aware that the expressions $A_0$, $A_1$ and $A_2$   provided above are those that appear in (2.10) of 
 \cite{Yuan2022}, but for computational purposes, we multiply by $-1$ and obtain the following  characteristic equation,   
 \begin{equation}\label{p_chin}
 \lambda^3 + A_2 \lambda^2 +A_1 \lambda +A_0 =0.
 \end{equation}
Keep in mind that this polynomial has similar roots and multiplicities as the one given in \eqref{p_mex}.

In \cite{Yuan2022}, the authors stated   that $A_2>0$, thus  the trace  $T(J(\overline{S}))$
 of $J(\overline{S})$ satisfies\footnote{If ${\mathcal A}$ is a $n\times n$ matrix, then  its characteristic polynomial is
$$p(\lambda)= \lambda^n - ({\rm tr} {\mathcal A}) \lambda^{n-1} + \cdots + (-1)^{n-1} ({\rm tr}({\rm adj} {\mathcal A})) \lambda
 + (-1)^n\det {\mathcal A},$$
 where ${\rm tr} {\mathcal A}$  and ${\rm adj} {\mathcal A}$ denote the trace and  the adjugate matrix of ${\mathcal A}$, respectively.} 
$$
T(J(\overline{S}))= -A_2 = -\alpha -\kappa_2 \overline{M}+ P(\overline{E}_2) U'_1(\overline{E}_2)<0.
$$
As it is well known, the eigenvalue sum is equal to the trace of the matrix, then $J(\overline{S})$ always 
has at least one eigenvalue that is negative.
Additionally, the determinant $D(J(\overline{S}) )$ of $J(\overline{S})$ is given by
$$
D(J(\overline{S}) )= -A_0 = -\alpha \kappa_2 \overline{M} \big[ P(\overline{E}_2) \big( U'_1(\overline{E}_2) -U'_2(\overline{E}_2)\big) +r_1\big],
$$
which is the product of the  three eigenvalues. Hence, if $J(\overline{S})$ has one zero eigenvalue,  then $D(J(\overline{S}) )=0$, leading to
$ U'_2(\overline{E}_2) -U'_1(\overline{E}_2)= \dfrac{r_1}{P(\overline{E}_2)}$.

\section{Bifurcation analysis}\label{sec_hopf}
The aim of the section is to discuss the conditions of the Hopf bifurcation in the model system \eqref{eq_model3}. Here, $\alpha$ is taken as the bifurcation parameter.

Let us  define $k_1=U'_1(\overline{E}_2)$, $k_2=U'_2(\overline{E}_2)$ and 
$\overline{k}= k_2-k_1= \dfrac{r_1}{P(\overline{E}_2)}$. 
The  points $\overline{S}=S_{12}$ and $\overline{S}=S_{23}$ nilpotent equilibria, since in these points $D(J(\overline{S}) )=0$ with
\begin{equation}\label{k2tilde}
\kappa_2=\widetilde{\kappa}_2 = \dfrac{acm^2\overline{E}_2}{(a+\overline{E}_2)^2 (a+ 2 \overline{E}_2) (r_1-\kappa_1 \overline{E}_2)}.
\end{equation}

\begin{remark}
It is worthwhile to point out that the authors in \cite{Yuan2022} misused the criteria of ${\mathbb R}^2$ for the existence of bifurcation points, 
being that   \eqref{eq_model3} is  a system in  ${\mathbb R}^3$.
\end{remark}
 The next proposition summarizes the local stability results and improves the flaws in item (3), Proposition 3.3 from \cite{Yuan2022}.
\begin{proposition}\label{prop_new}
The following statements on coexistence equilibria of \eqref{eq_model3}  hold, when exist.
\begin{enumerate}
\item  $S_1$ and $S_3$ are anti-saddle;
\item $S_2$ is a hyperbolic saddle;
\item if $J(S_{12 (23)})$ has only one eigenvalue zero and all of its eigenvalues are real, and  $T(J(S_{12 (23)}) )\neq 0$,  then $S_{12 (23)}$ is a saddle-node of codimension 1;
\item if $J(S_{1(3)})$  has one negative  eigenvalue and a pair of pure imaginary eigenvalues, then 
$T(J(S_{1(3)}) )< 0$ and   $D(J(S_{1(3)}) )<0$, so that $S_{1(3)}$ exhibits  a Hopf bifurcation point;
\item if $J(S_{12 (23)})$ has a double zero eigenvalue  and   $T(J(S_{12 (23)}) ) < 0$,  then at $S_{12 (23)}$ 
 occurs  a Bogdanov-Takens bifurcation.
\end{enumerate}
\end{proposition}

\subsection{Hopf bifurcation at $S_1$ or $S_3$}
To apply the Hopf bifurcation theorem, we must first determine the parameter values at which a pair of complex-conjugate eigenvalues of $J(S_{1(3)})$  cross, in a transversal manner, the imaginary axis of the complex plane while the other eigenvalue remains real and negative.

From Proposition \ref{prop_new}, we know that system \eqref{eq_model3} may exhibit Hopf bifurcation
at  $S_1$ or $S_3$, when exists, that is, if $T(J(S_{1(3)}) )=-A_2< 0$ and   $D(J(S_{1(3)}) )= -A_0 <0$.
These conditions mean that $A_2> 0$ and $A_0>0$.

For what follows in the rest of this paper, we will assume that the characteristic equation  \eqref{p_mex} of the system \eqref{eq_model3}
 has a pair of purely imaginary roots and  one real negative root.
Therefore,  the cubic equation   can be factored  as 
\begin{equation}\label{ch_eq1}
(\lambda + A_2)(\lambda^2+A_1)=0,
\end{equation}
imposing the usual condition for coefficients
\begin{equation}\label{AA3}
A_0=A_1A_2, \qquad A_1>0.
\end{equation}

Now, we choose the parameters $\alpha$ as the bifurcation parameter, and fix the other parameters at suitable values. 
By doing so, we define the  function
\begin{equation}\label{Ha}
H(\alpha)= A_0 - A_1 A_2 = h_2\alpha^2+h_1\alpha + h_0,
\end{equation}
where
\begin{align*}
h_0 &= - \frac{\delta_1 \overline{E}_2 + a r_2}{\kappa_2 (a+ \overline{E}_2)^5} \bigg[  2 \kappa_1  \overline{E}_2^2 + (a \kappa_1 + \delta_1 -r_1)\overline{E}_2 + a r_2\bigg] \\
 & \qquad \bigg[ 2 \kappa_1 \kappa_2 \overline{E}_2^4 + (5 a \kappa_1 \kappa_2 - r_1 \kappa_2) \overline{E}_2^3 +  2 a \kappa_2 (2 a \kappa_1- r_1) \overline{E}_2^2 + a \big( c m^2 - a r_1 + a^2 \kappa_1 \kappa_2 \big) \overline{E}_2  \bigg],\\
h_1 &=  \frac{1}{(a + \overline{E}_2)^2} \bigg[  -4 \kappa_1^2 \overline{E}_2^4 +  \big[ (6 r_1 - 4 \delta_1) \kappa_1 - 4 a \kappa_1^2 \big]\overline{E}_2^3 \\
& \qquad  -\big[ a^2 \kappa_1^2 +a \kappa_1 (2 c m-5 r_1+6 r_2)+2 r_1 (r_1 - \delta_1)+ \delta^2 \big]\overline{E}_2^2\\
& \qquad + \big[ (r_1 - 2 r_2) a^2 \kappa_1 - a(2 r_2 \delta_1 - 2 r_1 r_2 + r_1^2) \big] \overline{E}_2 - a^2 r_2^2  \bigg],\\
h_2 &= \frac{ -2 \kappa_1 \overline{E}_2^2 + (-a \kappa_1 - \delta_1 + 2 r_1) \overline{E}_2 + a(r_1 - r_2)}{a+ \overline{E}_2}.
\end{align*}
Clearly,  condition \eqref{AA3} is fulfilled when $H(\alpha)=0$. In fact, 
the discriminant of the quadratic polynomial in the variable $\alpha$ is given by
\begin{align*}
\Delta_h &= h^2_1 - 4 h_0 h_2 \\
 &= \frac{1}{(a+ \overline{E}_2)^4} \Big[ 4\kappa_1^2 \overline{E}_2^4 + \big(4 a\kappa_1^2 + 4 \delta_1 \kappa_1 - 6r_1 \kappa_1 \big) \overline{E}_2^3 \\
 &\qquad + \big[ \delta_1^2 + 2a \kappa_1 \delta_1 - 2 r_1 \delta_1 + (r_1 - a\kappa_1)^2 -a (3 r_1 -4r_2)\kappa_1 + r_1^2 \big] \overline{E}_2^2 \\
&\qquad + a \big[ r_1^2 - 2r_1 r_2 + 2r_2 \delta_1 - ar_1 \kappa_1 + 2a r_2 \kappa_1 \big] \overline{E}_2 + a^2 r_2^2 \Big)\\
&\qquad - \frac{4 \overline{E}_2 (\delta_1 \overline{E}_2 + ar_2)}{(a+ \overline{E}_2)^6 \kappa_2} 
 \Big[ \big(2 \kappa_1 \overline{E}_2 + a \kappa_1 + \delta_1 - 2r_1 \big) \overline{E}_2 +  a( r_2 - r_1 ) \Big] \\
&\qquad   \Big[ \big( 2 \kappa_1 \overline{E}_2 + a \kappa_1 + \delta_1 - r_1 \big)\overline{E}_2 + a r_2 \Big]
 \Big[ (\overline{E}_2 + a)^2 \big( (a + 2 \overline{E}_2)\kappa_1- r_1\big) \kappa_2+ a c m^2 \Big].
\end{align*}
By using the quadratic formula, we can see that the roots of $H(\alpha)=0$ yield
\begin{equation*}
\alpha^- = \frac{-h_1 - \sqrt{\Delta_h}}{2 h_2}, \qquad \alpha^+ = \frac{-h_1 + \sqrt{\Delta_h}}{2 h_2}.
\end{equation*}

Now, we remark that  if condition $\Delta_h=0$ holds, then there is only one real root for $H(\alpha)=0$. This reads
\begin{align*}
\alpha^* =- \frac{h_1}{2 h_2} 
=&  \frac{1}{2 (a + \overline{E}_2) \big( \overline{E}_2\kappa_1(a + 2 \overline{E}_2)+ (\delta_1 - 2r_1)\overline{E}_2 +ar_2 - ar_1 \big)} \\
& \qquad \Big[ - 4\kappa_1^2 \overline{E}_2^4
 - 2\kappa_1 \big( 2a\kappa_1+ \delta_1 - 3r_1 \big) \overline{E}_2^3 \\
& \qquad - \big( a^2 \kappa_1^2 + a(2 cm - 5r_1 + 6r_2) +2 r_1 (\delta_1 + r_1) + \delta_1^2 \big) \overline{E}_2^2 \\
& \qquad + a \big( -r_1^2 + r_1(2r_2 + a\kappa_1)- 2r_2 (a\kappa_1 + \delta_1) \big) \overline{E}_2 - a^2 r_2^2
)\Big].
\end{align*}

Let us turn our attention to determine whether $A_1>0$ occurs.
We calculated this by using Mathematica software and found that this happens when
$$r_1 >\dfrac{(a+E_2)(cm E_2+(a+E_2)r_2)\alpha}{cm E_2^2 + (a+E_2)(a\alpha + E_2(r_2+2\alpha))},$$
$$
0 < \kappa_1 < \dfrac{cmE_2\big[E_2(r_1-\alpha)-a\alpha\big] + (a+E_2)\big[a(r_1-r_2)\alpha -r_2\alpha E_2+E_2 r_1(r_2+2\alpha)\big] }{E_2(a+2E_2)\big [ c m E_2 + (a +E_2) (r_2+\alpha)\big]},
$$
and $\widehat{\kappa}_2  < \kappa_2$, with
$$
\widehat{\kappa}_2 = \frac{ - a c m^2  \overline{E}_2 (\delta_1 \overline{E}_2 + a r_2)}{(a + \overline{E}_2)^2 (\widetilde{q}_3 \overline{E}_2^3 + \widetilde{q}_2 \overline{E}_2^2 + \widetilde{q}_1 \overline{E}_2 + \widetilde{q}_0)},
$$
\begin{align*}
\widetilde{q}_3 &=  2\kappa_1 (\delta_1 + \alpha), \\
\widetilde{q}_2 &= a \kappa_1 \big(cm + 3\alpha +3r_2 \big)- 2\alpha r_1 + \delta_1(\alpha-r_1), \\
\widetilde{q}_1 &= a \big( a \kappa_1 (r_2 + \alpha)+ (cm - 3r_1 +2r_2) \alpha -r_1 r_2 \big),\\
\widetilde{q}_0 &=  a^2 \alpha (r_2 - r_1).
\end{align*}

\begin{remark}
We would like to point out that after double-checking the results, we discovered that  $\widetilde{q}_2$, $\widetilde{q}_1$ and $\widetilde{q}_0$ expressions  obtained are different from those in \cite{Yuan2022} on page 112. 
Also,  in Proposition 3.3 of \cite{Yuan2022}, the authors failed to make the assumption that $\kappa_2 < \widehat{\kappa}_2$. 
This claim will be supported in Section \ref{sec_nume}, where some numerical simulations of the system \eqref{eq_model3} for various parameter values will demonstrate the formation of limit cycles via Hopf bifurcations in equilibrium, but using the proper inequality
$\widehat{\kappa}_2 < \kappa_2 $. 
\end{remark}

Let us assume that $\hat{\lambda}_1$, $\hat{\lambda}_2$, $\hat{\lambda}_3$ are the eigenvalues of \eqref{matrizJ}, 
 or equivalently,  the roots of equation \eqref{ch_eq1}. Thus, there exists an eigenvalue of 
\eqref{p_mex}, say $\hat{\lambda}_1$,  such that  $\hat{\lambda}_3=\overline{\hat{\lambda}}_1$, 
and the other eigenvalue  satisfies  $\hat{\lambda}_2<0$. Consequently, $\hat{\lambda}_2 + \hat{\lambda}_3 = -A_2-\hat{\lambda}_1$.
As  $\hat{\lambda}_1 $ and $\hat{\lambda}_3$ are complex conjugates, it follows that  $2 \gamma = -A_2-\hat{\lambda}_1$,
where $\gamma $ is the real part of  $\hat{\lambda}_1$ and $\hat{\lambda}_3$.  

To  assure the occurrence of the Hopf bifurcation we need  to guarantee
the transversality condition of the Hopf bifurcation theorem.
By using \eqref{AA3} we have $2 \gamma = -\dfrac{A_0}{A_1} - \hat{\lambda}_1$, and so 
 $\gamma = -\dfrac{1}{2}\left(\frac{A_0}{A_1} + \hat{\lambda}_1\right)$.
Differentiating this equation  with respect to $\alpha$,  and arranging the terms, one reaches
{\scriptsize $$
 \dfrac{\partial\gamma}{\partial \alpha}\bigg|_{\alpha=\alpha^\pm}  =
\dfrac{\overline{E}_2 (\delta_1 \overline{E}_2 + r_2a)^2\Big[\kappa_2 (a+\overline{E}_2)^2(-r_1 +(a+2 \overline{E}_2)\kappa_1+acm^2)\Big]\Big[ (a + \overline{E}_2)^2 (a + 2 \overline{E}_2) (-r_1 + \overline{E}_2 \kappa_1) \kappa_2+a c m^2\overline{E}_2  \Big]}{\Big[(\overline{q}_3\overline{E}_2^3   + \overline{q}_2\overline{E}_2^2+\overline{q}_1\overline{E}_2+\overline{q}_0)(a+\overline{E}_2)^2\kappa_2+ acm^2(\delta_1 \overline{E}_2+r_2a) \overline{E}_2\Big]^2}
$$}
where 
\begin{align*}
&\overline{q}_3 = 2 (\alpha^\pm+\delta_1)\kappa_1,\\
&\overline{q}_2 = (3a\kappa_1+\delta_1)\alpha^\pm + a\kappa_1 (3 r_2 + cm)-r_1\delta_1,\\ 
&\overline{q}_1 = a(cm + r_1+2r_2+a\kappa_1)\alpha^\pm + a r_2(-r_1+a\kappa_1),\\
&\overline{q}_0 = a^2(r_2+r_1)\alpha^\pm.
\end{align*}
After some algebraic computations, we get that
$\dfrac{\partial\gamma}{\partial \alpha}\bigg|_{\alpha=\alpha^\pm} =0$ if and only if 
\begin{equation}\label{hh1}
a c m^2 \overline{E}_2+(a+\overline{E}_2)^2 (a+2\overline{E}_2) (-r_1+\kappa_1 \overline{E}_2)\kappa_2=0,
\end{equation}
or
\begin{equation}\label{hh2}
a c m^2+(a+\overline{E}_2)^2 (-r_1+(a+2\overline{E}_2) \kappa_1)\kappa_2 =0.
\end{equation}

Through an immediate calculation,  from \eqref{hh1} we obtain  that
$$
\kappa_2 = \dfrac{acm^2\overline{E}_2}{(a+\overline{E}_2)^2 (a+ 2 \overline{E}_2) (r_1-\kappa_1 \overline{E}_2)},
$$
that is,   $\kappa_2=\widetilde{\kappa}_2$   where   $D(J(\overline{S}) )=0$ for $\overline{S}= S_{12}$ or  $\overline{S}= S_{23}$.
On the other hand, the equation \eqref{hh2} is true for 
\begin{equation}\label{k2s123} 
\kappa_2 =\dfrac{acm^2}{(a+\overline{E}_2)^2(r_1-a\kappa_1-2 \kappa_1 \overline{E}_2)}.
\end{equation}

A straightforward computation shows that the equilibrium $\overline{S}=S_{123}$  exists if equalities
$U(E_2)=0$, $\dfrac{dU}{dE_2}=0$ and $\dfrac{d^2U}{dE_2^2}=0$ are fulfilled simultaneously, yielding
$$
E_2^*=a\left( -1 +\dfrac{3cm}{\delta_1}\right), \qquad  \kappa_1 = \kappa_1^*=\dfrac{r_1\delta_1}{a\delta_2}, \qquad 
\kappa_2 = \kappa_2^*= \dfrac{\delta_1^2\delta_2}{27ac^2mr_1}.
$$
Now, after substituting  $\overline{E}_2=E_2^*$ and  $\kappa_1 = \kappa_1^*$ into \eqref{k2s123} we get  that $\kappa_2 = \kappa_2^*$.
So,   \eqref{k2s123}  is associated with  $\overline{S}= S_{123}$.
Therefore, $\dfrac{\partial\gamma}{\partial \alpha}\Big|_{\alpha=\alpha^\pm}\neq 0$ and
 the condition of transversality is verified. Hence, the system  \eqref{eq_model3}  exhibits  a Hopf bifurcation at $\alpha^\pm$.

\begin{remark}
We observe that  the  factor $$a c m^2 \overline{E}_2+(a+\overline{E}_2)^2 (a+2\overline{E}_2) (-r_1+\kappa_1 \overline{E}_2)\kappa_2$$
 is missing in the algebraic expression of  $\dfrac{\partial\gamma}{\partial \alpha}\Big|_{\alpha=\alpha^\pm}$  given in \cite{Yuan2022}. Nevertheless, the authors  claim that $\dfrac{\partial\gamma}{\partial \alpha}\Big|_{\alpha=\alpha^\pm}=0$ if only if 
$$ r_1=\kappa_1 (a+ 2 \overline{E}_2)+\dfrac{acm^2}{(a+\overline{E}_2)^2\kappa_2},$$
equivalently to \eqref{hh2},  occurs only for points $\overline{S}= S_{12}$ or  $\overline{S}= S_{23}$ or  $\overline{S}= S_{123}$. 
The authors are incorrect in their assertion because this value of $r_1$ only appears in $\overline{S}= S_{123}$ and not in $\overline{S}= S_{12}$ or  $\overline{S}= S_{23}$.
\end{remark}

As a result, we have proved the following theorem, which outlines the conditions that are necessary for a Hopf bifurcation to occur.
\begin{theorem}\label{teo1}
For the set of differential equations \eqref{eq_model}, if $\widehat{\kappa}_2  < \kappa_2$ and  the coexistence point $S_1$ or $S_3$ exists, then
\begin{enumerate}
\item[{1.}] if $\alpha=\alpha^-$, or $\alpha=\alpha^+$,   system \eqref{eq_model} goes through  a Hopf
bifurcation;
\item[{2.}] if $\alpha=\alpha^*$,   system \eqref{eq_model} goes through a Hopf bifurcation, which can be degenerate.
\end{enumerate}
\end{theorem}
In the next section we provide numerical examples to highlight the existing Hopf bifurcation 
of codimension 1 and 2.

\section{Numerical experiments}\label{sec_nume}
In this section we illustrate the creation of  limit cycles emanating from the Hopf bifurcations at the equilibrium.

We would like to point out that the reader should be aware that the main purpose of the examples provided is not to discover some new interesting bifurcation behaviors, but rather to demonstrate the validity of our version of Theorem \ref{teo1}, which corrects what is stablished in   Proposition 3.3 in \cite{Yuan2022}.

The bifurcation analysis is done using MatCont \cite{matcont}. This is a package for numerical continuation that runs in a Matlab environment.
This numerical simulation tool  can also be used to determine the first ($\ell_1$) and second ($\ell_2$) Lyapunov coefficients by applying the formulae proposed by Kuznetsov in  \cite{Kuznetsov_libro}.

For the purpose of this section, we consider a case that looks slightly simpler case than the general case, where
$ a=c=m=1$,  which implies that systems \eqref{eq_model} and \eqref{eq_model2} are the same. %
Biologically this means that the number of removal of adult individuals born from {\em E. onukii} attacked by one {\em A. baccarum} is equal to $\dfrac{E_2M}{1+E_2}$.  Notwithstanding the system still exhibits very complex dynamics even under this restriction.

\begin{figure}[h!]
\centering
\includegraphics[scale=0.36]{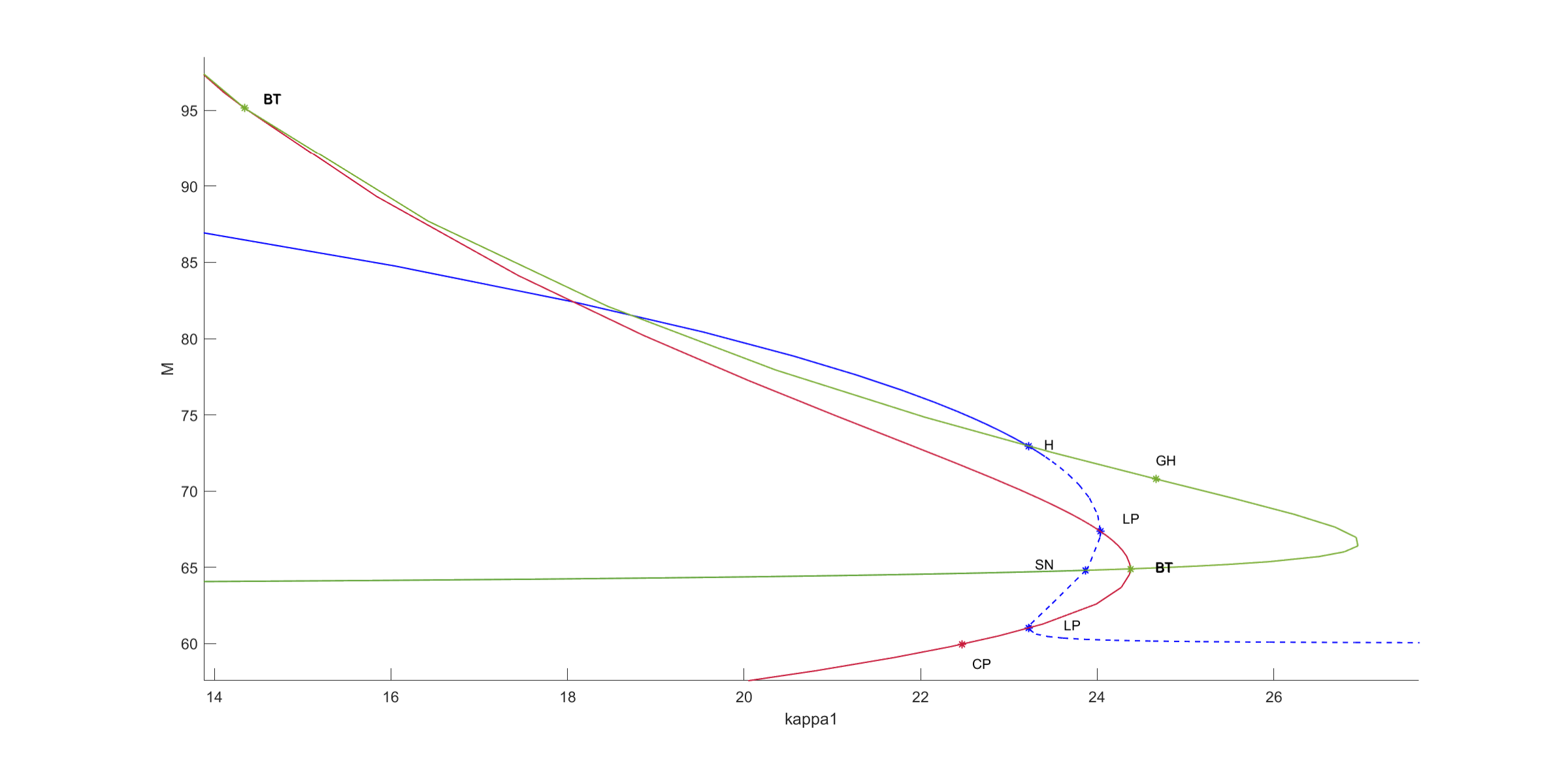} 
\caption{This plot shows the continuation of the equilibrium curve detected with MatCont in  the $(\kappa_1,M)$-plane where the parameter values $r_1=1.6$, $r_2=60$, $\alpha=51.57$  and $\kappa_2=0.026671$ stand fixed:
switching from one branch of equilibria to the other and finding two branch points on the other branch.
BT, LP, SN, H, CP  and GH denote  Bogdanov-Takens bifurcation point,  limit point, saddle-node point, Hopf bifurcation point, cusp point and degenerate Hopf  bifurcation point, respectively. The blue solid and dashed lines represent the stable and unstable hopf bifurcation, respectively.}\label{fig_bifg}
\end{figure}

In Figure \ref{fig_bifg} we see  the continuation of the equilibrium curve detected with MatCont in  the $(\kappa_1,M)$-plane.
During the bifurcation analysis some parameters are fixed at $r_1=1.6$, $r_2=60$, $\alpha=51.57$ and $\kappa_2=0.026927991$ because for this values most types of activity can be seen.
The branch of Hopf bifurcations is continued in two parameters $(\kappa_1,\kappa_2)$, and a generalized Hopf (GH) point is located.
Switching from one branch of equilibria to another branch of equilibria.

And then, one has to note that the  bifurcations and related dynamics, including Hopf (H and GH) and Bogdanov-Takens (BT) bifurcations, which reveals the complex dynamics behaviors and the reason behind the complexity of the model \eqref{eq_model}.

\subsection{Hopf bifurcation}
\begin{itemize}
\item The system \eqref{eq_model} goes through a Hopf bifurcation at the equilibrium point $\overline{S}$, whose coordinates are
$$(\overline{E}_1,\overline{E}_2,\overline{M})=(  0.614426554662767,  0.528099623732648, 72.9478611523887),$$
and parameter values  $\alpha = 51.57$, $\kappa_1 = 23.2197961461739$, $\kappa_2= 0.026671$ 
$r_1=60$ and  $r_2=1.6$.
The  eigenvalues of the Jacobian matrix  \eqref{matrizJ} evaluated at   $\overline{S}$ are given by 
$$ \lambda_1= -109.28094024998, \quad  \lambda_2=1.09824705073231i, \quad \lambda_3=-1.09824705073231i .$$ 
Moreover, 
$\widehat{\kappa}_2=0.0000801963 < \kappa_2=0.026671$, so  the hypothesis of Theorem \ref{teo1} holds.

The  first Lyapunov coefficient, obtained by using the software MatCont corresponds to  $\ell_1= 0.02036690 >0$. Its positivy means that the Hopf bifurcation is subcritical, namely, the limit cycle arising near equilibrium point is unstable. 
As indicated in Figure \ref{hopp1}, the equilibrium point is stable for $\kappa_1  < 23.2197961461739$, afterwards loses stability at the Hopf (H) bifurcation point, and then degenerates into instability. 

We plot the phase portrait in Figure \ref{hopp2}, where the subcritical Hopf bifurcation takes place, an unstable limit cycle  (LPC, colored red) mount from the unstable equilibrium.
\begin{figure}[hbt]
\includegraphics[scale=0.36]{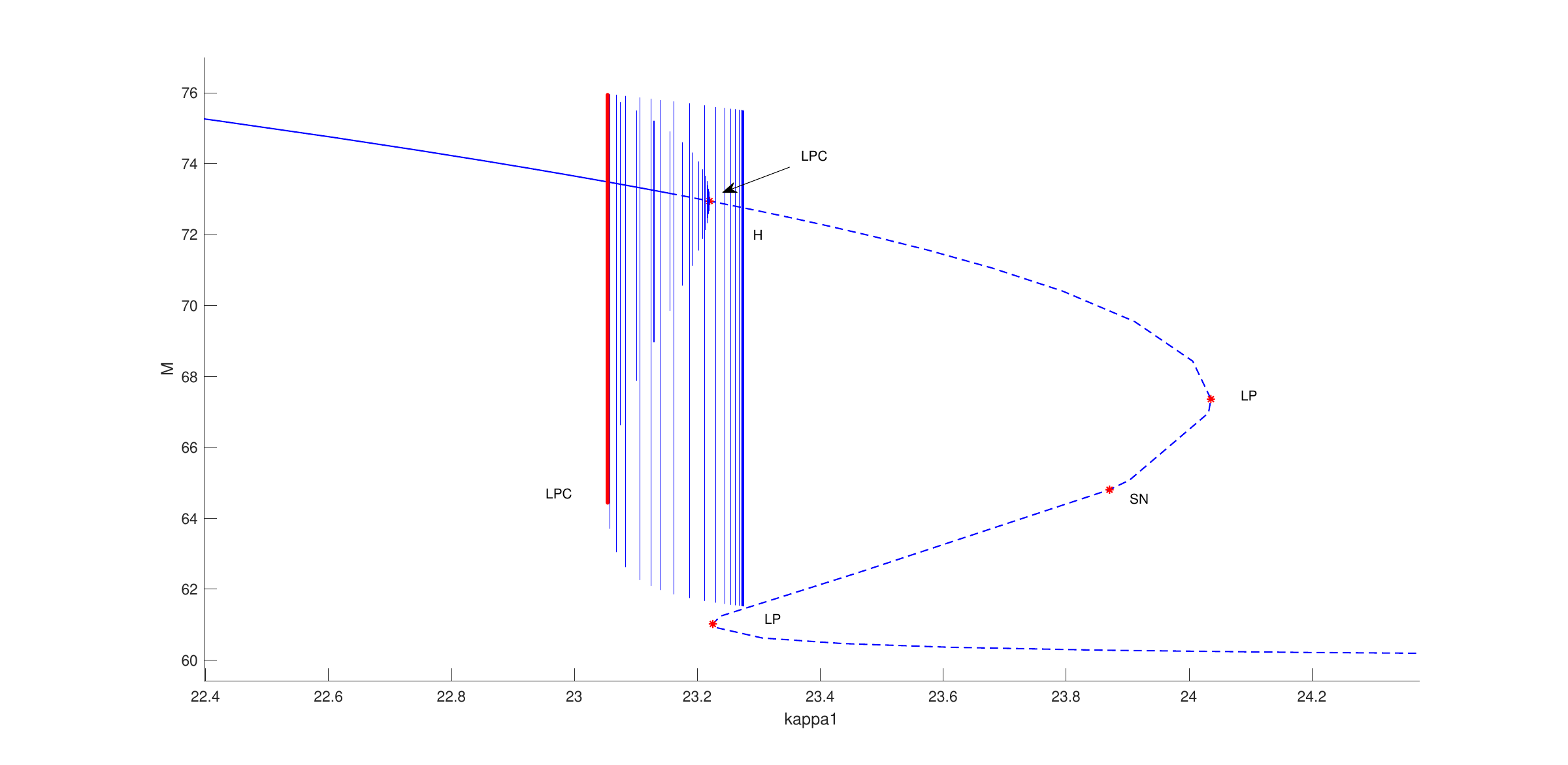}
\caption{Diagram of  numerical continuation of  subcritical Hopf (H)  point in the $(\kappa_1,M)$-plane. 
The stability  of the equilibrium point changes from the stable (solid line)  to the unstable (dashed line). This 
behavior corresponds to the subcritical Hopf bifurcation.}
\label{hopp1}
\end{figure}

\begin{figure}[hbt]
\includegraphics[scale=0.35]{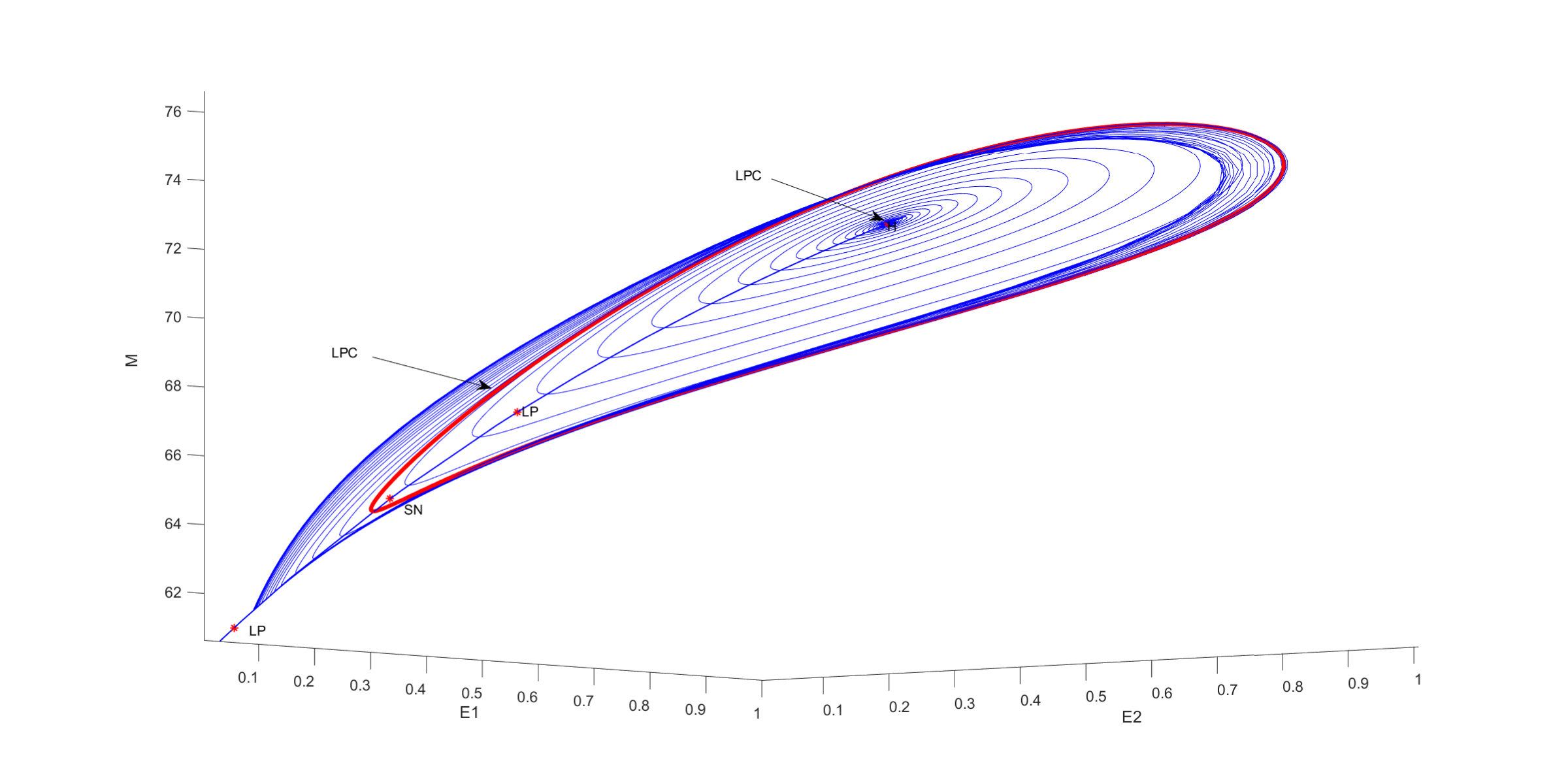}
\caption{$(E_1,E_2,M)$-phase space. Shown an unstable limit cycle  (LPC) solution of system \eqref{eq_model} obtained using MatCont by continuation of the orbit emerging  from subcritical Hopf (H) point
$(\overline{E}_1,\overline{E}_2,\overline{M})=(  0.614426554662767,  0.528099623732648, 72.9478611523887)$.}\label{hopp2}
\end{figure}
\end{itemize}

\begin{itemize}
\item The model \eqref{eq_model} undergoes a subcritical Hopf bifurcation at
$$(\overline{E}_1,\overline{E}_2,\overline{M})=(0.702380773962441,0.418254438462965,  73.7892484519822)$$
for $m=1$, $a=1$, $r_1=62.27545$, $r_2=1.6$,  $\kappa_1= 24.5$,  $\kappa_2=0.02568$, $\alpha= 37.083850149878$,

Since $\widehat{\kappa}_2=  0.000113965< \kappa_2=0.02568$, as a result, Theorem \ref{teo1}  is true. Given that
$\ell_1 = 0.004838562$, this implies that the bifurcating limit cycle is  unstable. 
Figure \ref{hopp22} displays  the phase portrait for this situation.
\begin{figure}[hbt]
\includegraphics[scale=0.35]{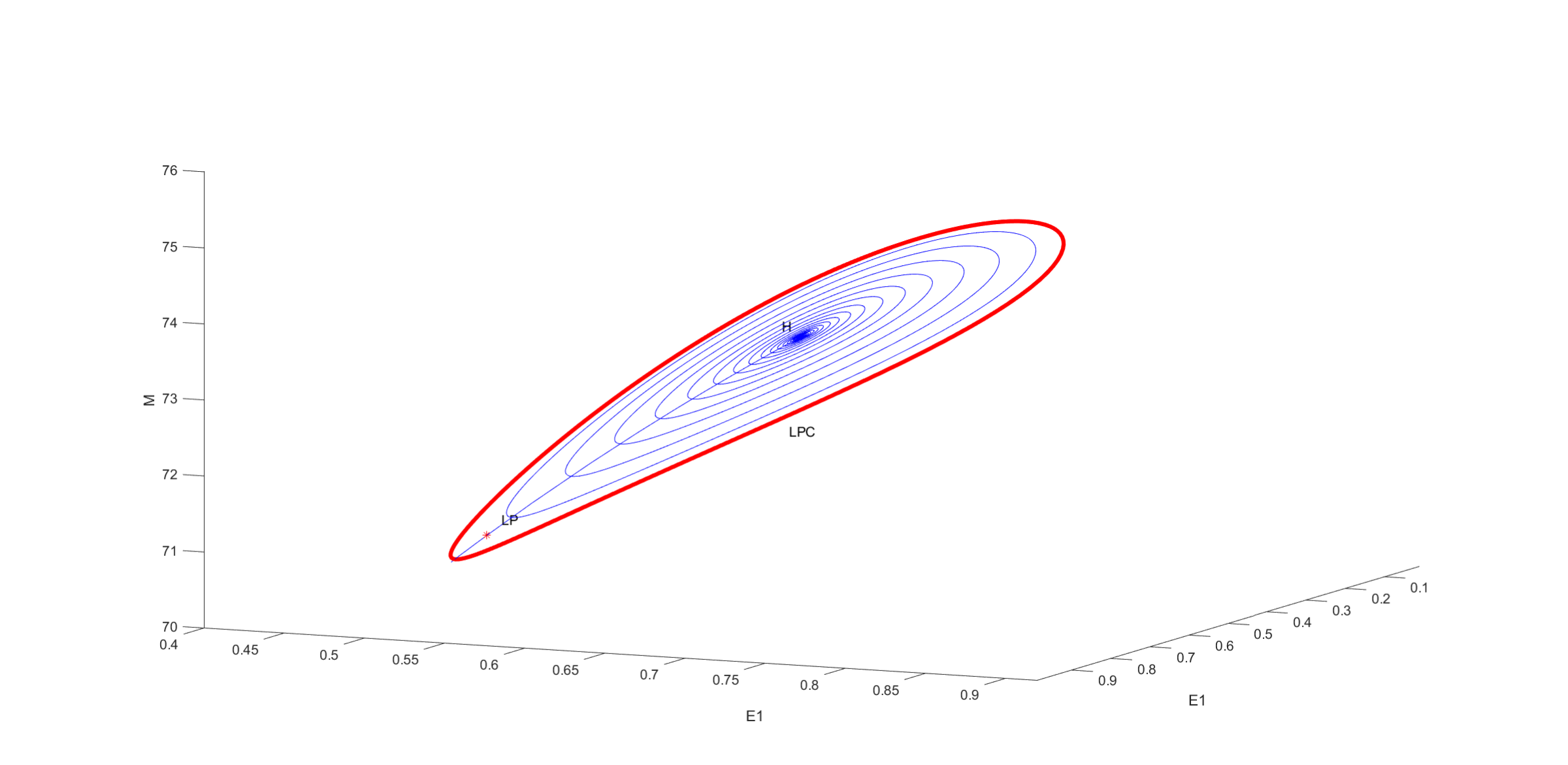}
\caption{In the $(E_1,E_2,M)$-phase space, it is shown an unstable limit cycle  (LPC) solution of the system 
emerging  from subcritical Hopf (H) point
$(\overline{E}_1,\overline{E}_2,\overline{M})=(0.702380773962441,0.418254438462965,  73.7892484519822)$.}\label{hopp22}
\end{figure}
\end{itemize}

\subsection{Bautin bifurcation}
An important factor in understanding a dynamical system's overall behavior is the presence of a codimension-2 bifurcation point, which has a significant impact on the qualitative behavior of the system.
Here, we perform a bifurcation analysis of Bautin or degenerate  codimension-2 Hopf bifurcation
at which a fold limit cycle and a Hopf bifurcations occur together.

\begin{figure}[hbt]
\centering
\includegraphics[scale=0.5]{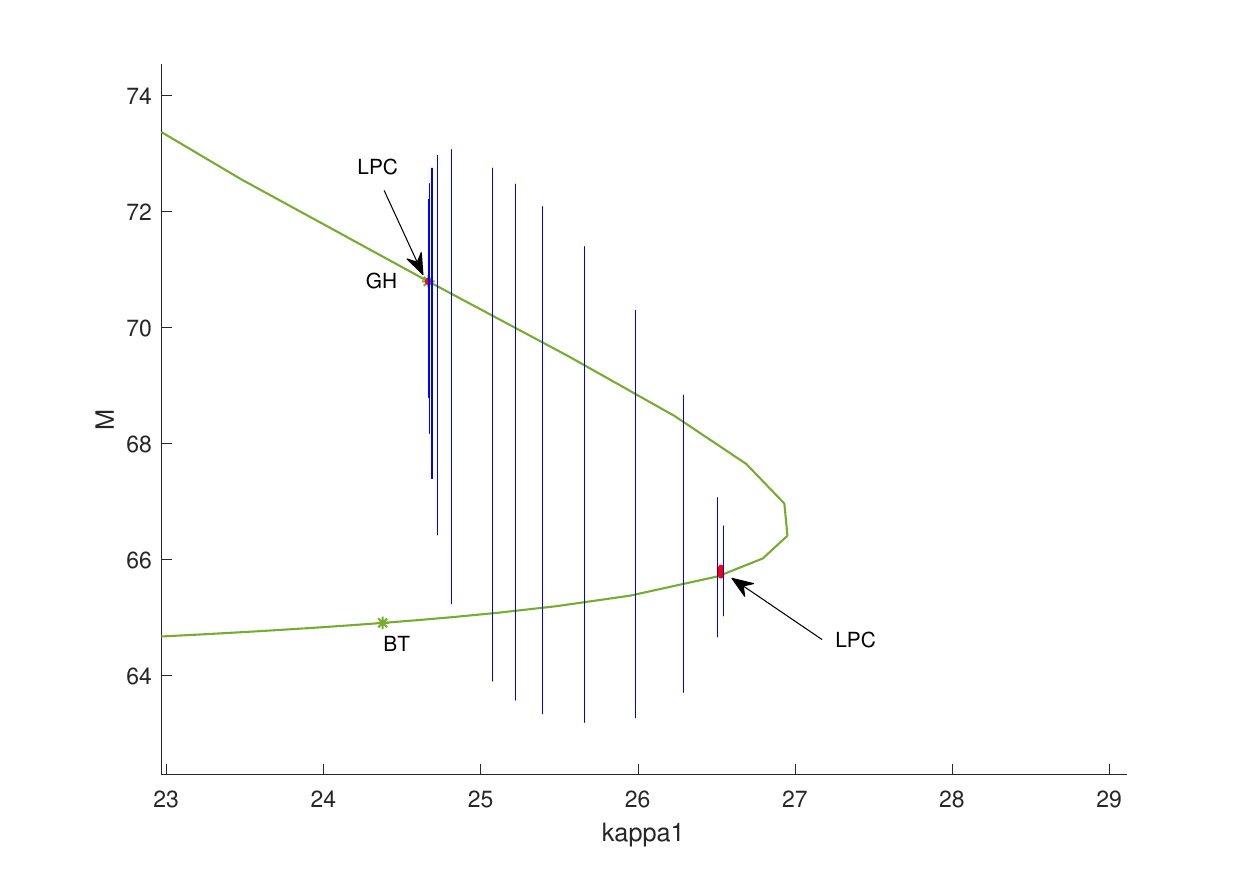}
\caption{Diagram of  numerical continuation of generalized Hopf (GH)  point in the $(\kappa_1,M)$-plane showing the  periodic orbits that emerge from such  point,
when $r_1=60$, $r_2=1.6$, $\alpha = 51.57$ $\kappa_2=0.026671$.
 BT is Bogdanov-Takens point and  LPC is limit point of cycle.}\label{fig_bautin1a} 
 \end{figure}

The  numerical simulation tool MatCont indicates that a Bautin point takes place at
$$(\alpha,\kappa_1,\kappa_2,r_1,r_2)=(53.1351, 24.665343,  0.026927991, 60,  1.6).$$
As required by the Theorem \ref{teo1}, it is true $\widehat{\kappa}_2= 0.0000847597 < \kappa_2 =0.026927991$.
For the chosen parameters in Figure \ref{fig_bautin1a}, a 2D figure in the $(\kappa_1,M)$-plane displaying 
 the  periodic orbits that emerge from the generalized Hopf.

\begin{figure}[hbt]
\subfigure[Phase portrait exhibiting Bautin point of codimension 2.]{\includegraphics[scale=0.35]{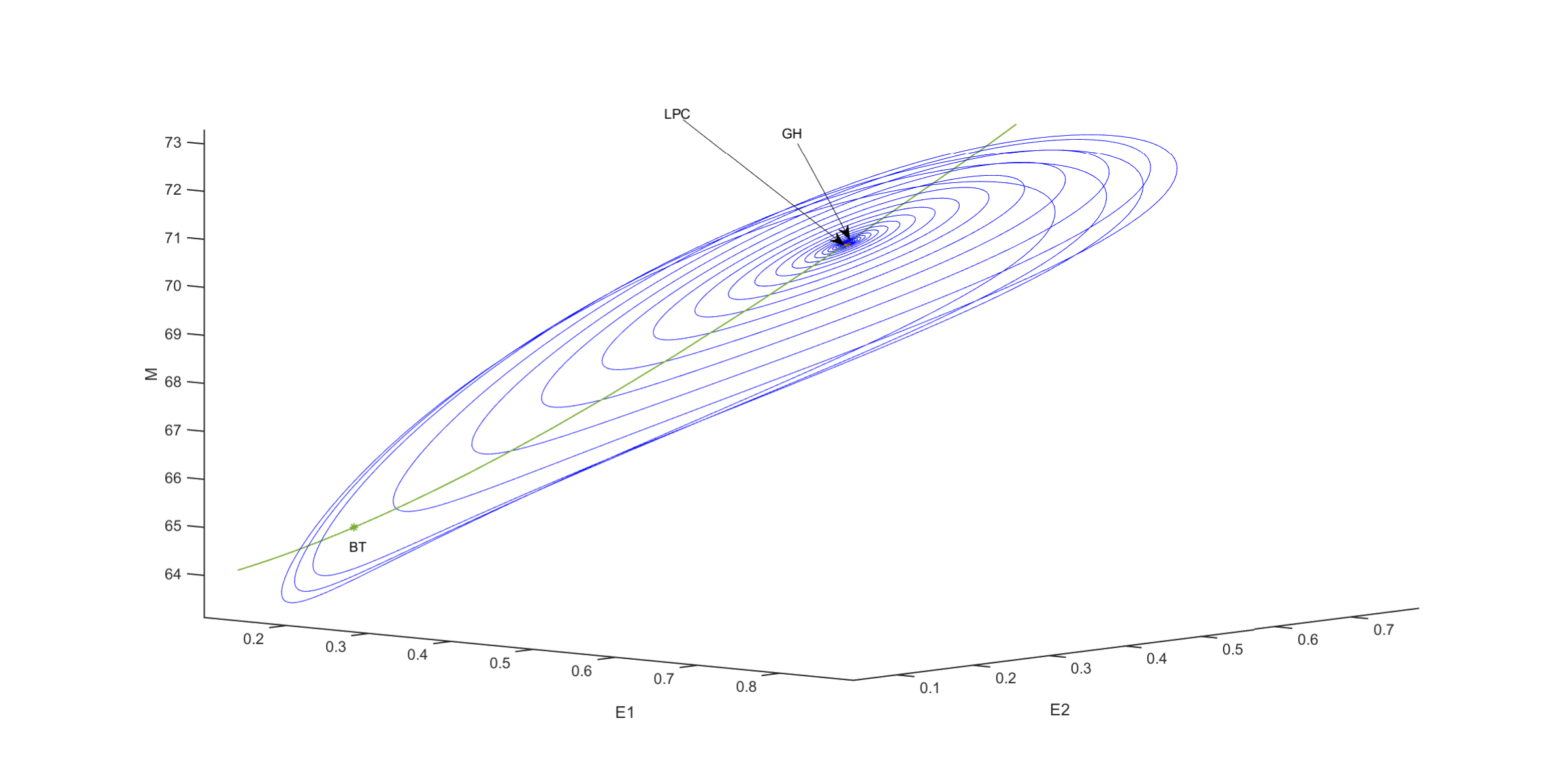}}
\subfigure[Small limit cycle bifurcating from  the generalized Hopf point.
 Zoom plot of  (a).]{\includegraphics[scale=0.3]{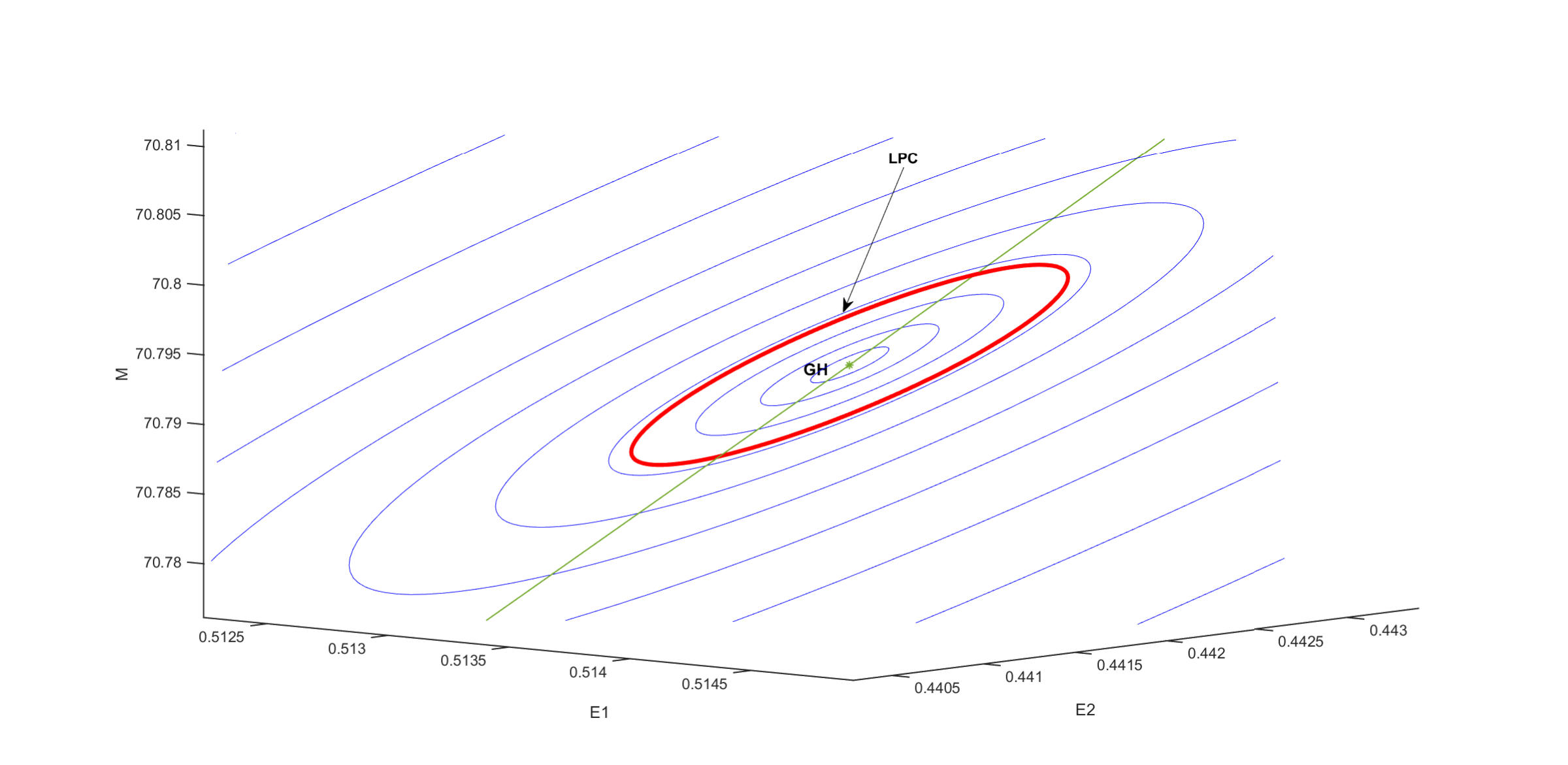}	}
\caption{$(E_1,E_2,M)$-phase space. The generalized Hopf bifurcation point  is indicated by GH, while  the 
limit point of cycle   (labeled  by LPC)   is a fold bifurcation of the cycle.}\label{fig_bautin1}
\end{figure}

The first Lyapunov coefficient $\ell_1$ shrinks to zero, which means that the Hopf bifurcation of 
 $$(\overline{E}_1,\overline{E}_2,\overline{M})=(0.51386752, 0.44166914,  70.794719)$$ 
degenerates \cite{perko}, at which, the eigenvalues are given by
$$\lambda_1=-109.326162726243, \quad \lambda_2= 1.0910208605339 i, \quad  \lambda_3= -1.0910208605339 i.$$
In order to see if this  point  is a Bautin bifurcation  or has a higher order degeneracy, 
we need to calculate the second Lyapunov coefficient $\ell_2$
with $\kappa_1$ and  $\kappa_2$  as unfolding parameters.
With the help of MatCont we obtain that  $\ell_2= -0.001914430  \neq 0$. 
In Figure \ref{fig_bautin1} the result of the bifurcation analysis is shown. Notice that 
LPC point (which correspond to the collapse of stable-unstable limit cycles) and 
Bautin (generalized Hopf,  GH) bifurcation point   are labeled.

\subsection{Degenerate Hopf bifurcation  with four  limit cycles}
We provide one example to illustrate the existence of four limit cycles that bifurcate from a degenerate Hopf bifurcation  for parameter values  close to the focus equilibrium point.

The system \eqref{eq_model} has an equilibrium at
$$(\overline{E}_1, \overline{E}_2, \overline{M})= (0.667957076496141, 0.37421233014172, 72.9092840794877)$$
corresponding to
 $m=1$, $a=1$, $r_1=62.27545$, $r_2=1.6$,  $\kappa_1= 24.638748097512 $,  $\kappa_2=0.02568$,  and $\alpha= 34.888830547725$,
with eigenvalues 
$$
\lambda_1= -93.80919627015, \quad  \lambda_2=0.519624557163896i,  \quad  \lambda_3=-0.519624557163896i.$$
Since
$\widehat{\kappa}_2= 0.000117406  < \kappa_2 = 0.02568$, then the numerical results agree quite well with the analytical ones
for the existence of a Hopf bifurcation.

\begin{figure}[h!]
\centering
\includegraphics[scale=0.35]{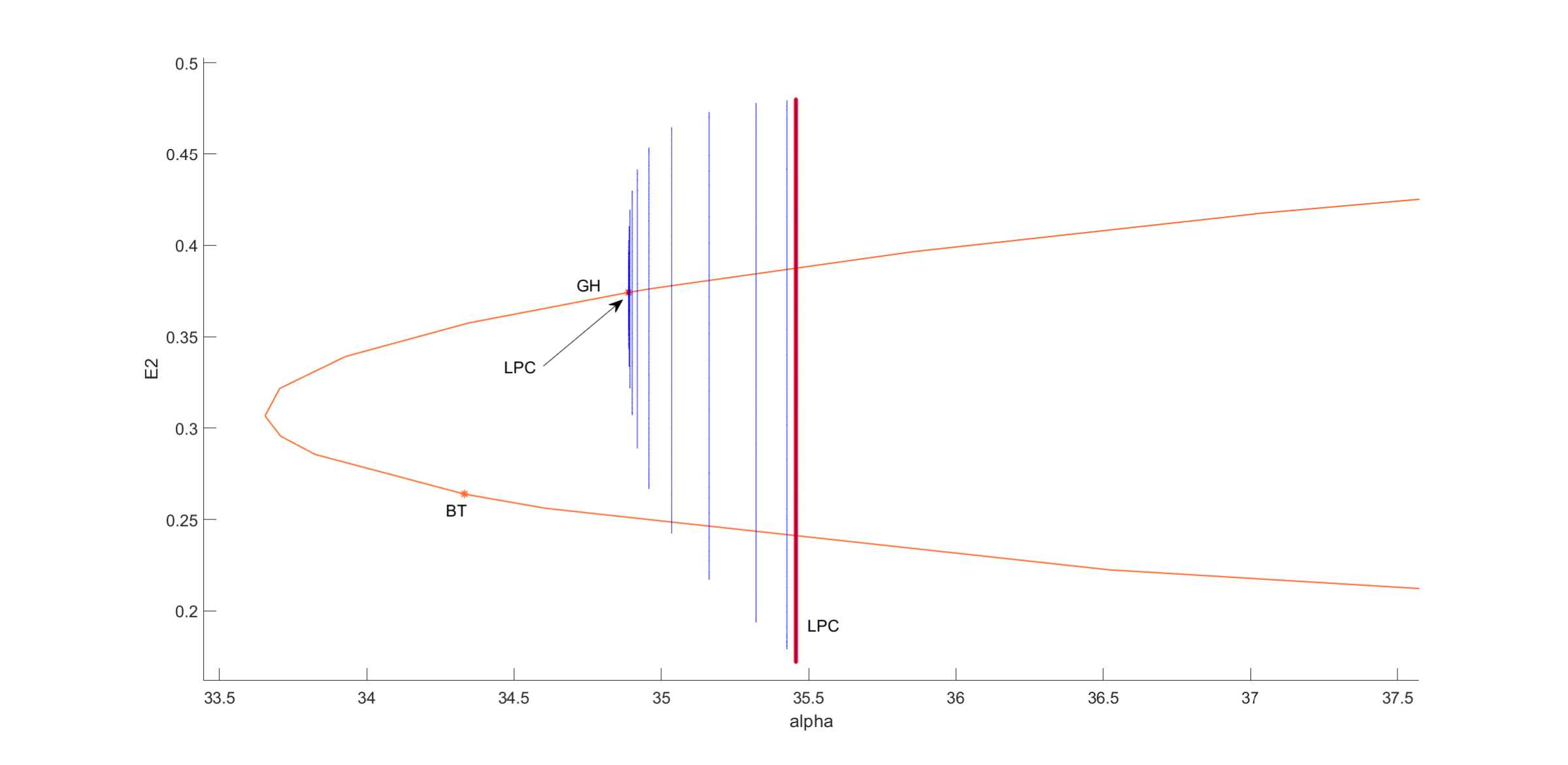} 
\caption{Bifurcation diagram in  $(\alpha, E_2)$-plane with $r_1=62.27545$, $r_2=1.6$,  $\kappa_1= 24.638748097512 $,  $\kappa_2=0.02568$.}\label{bautdeg1}
\end{figure}

The numerical results provided by MatCont says that the corresponding first Lyapunov coefficients $\ell_1$ becomes zero. This implies that the Hopf bifurcation is degenerate \cite{perko}, but  the critical second Lyapunov coefficient is different from zero, namely $\ell_2 = 0.001968325$. By choosing the parameters $\alpha$ and  $\kappa_1$ the fold limit cycle codimension-2 Hopf bifurcation occurs at  
$$(\alpha, \kappa_1)= (34.332136490836, 24.806008381396).$$
The projection of Hopf bifurcation regions on the $(E_2,\alpha)$-plane is shown in Figure \ref{bautdeg1}, in which $\ell_1=0$  and $\ell_2$ are the projection of Bautin bifurcation lines. The bifurcation portrait and the
corresponding phase portraits is shown in Figure \ref{bautdeg1},  where it is seen  that
the degeneracy gives rise to limit folding cycles (LPCs) that surround the equilibrium point.

\begin{figure}[hbt]
\subfigure[Phase portrait exhibiting Bautin point of codimension 2 for parameter values close to the  nilpotent singularity's focus  of codimension 3.]{\includegraphics[scale=0.35]{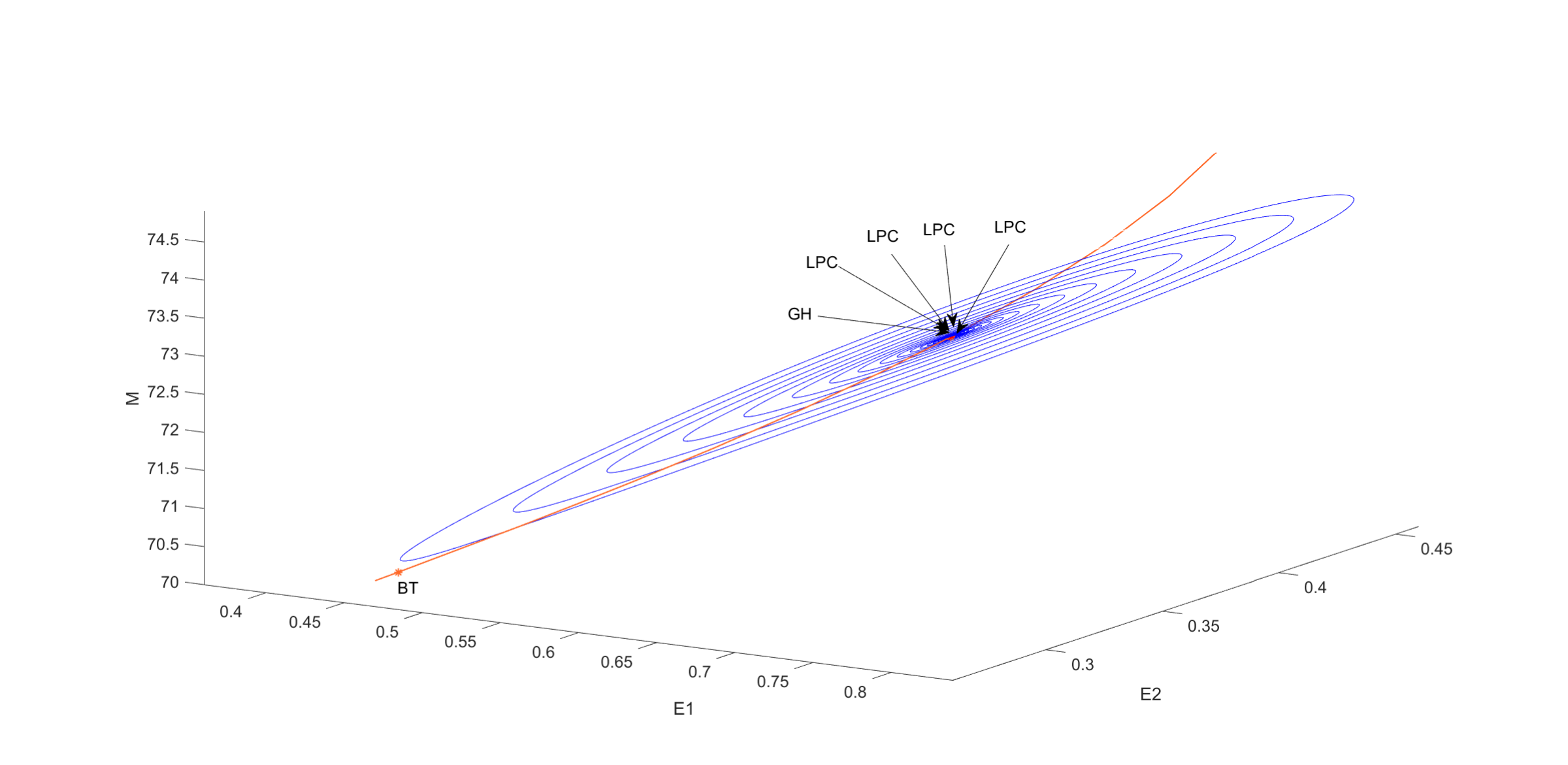}}
\subfigure[Phase portraits near the Hopf bifurcation of codimension 2. Four small limit cycle bifurcating from  the Hopf point. Zoom plot of  (a).]{\includegraphics[scale=0.3]{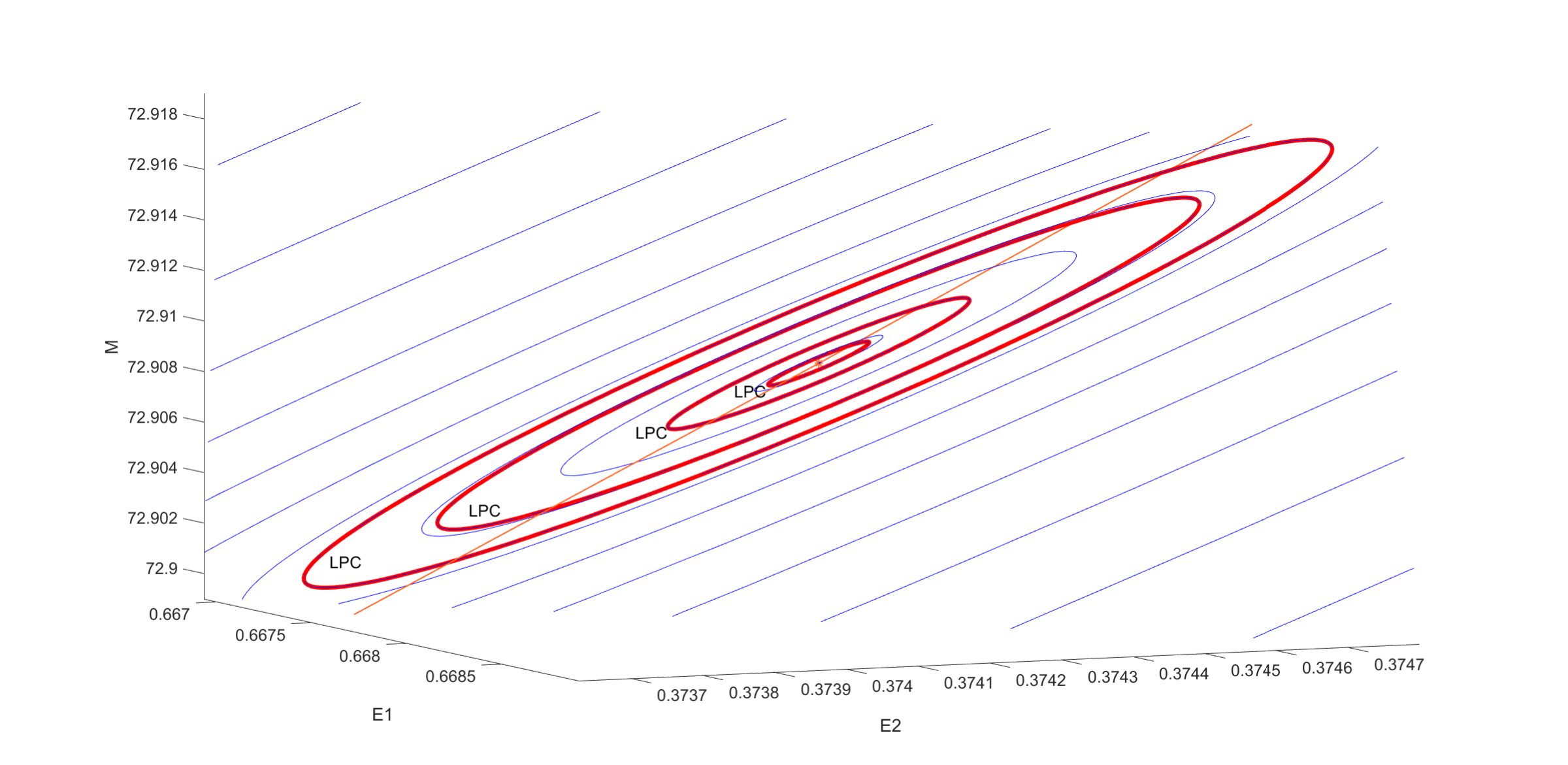}}
\caption{In the $(E_1,E_2,M)$-phase space. We show a codimension 2-Hopf  bifurcation of system \eqref{eq_model}. The generalized Hopf bifurcation point  is indicated by GH and   the small four
 limit point of cycle  are labeled  by LPC.}\label{bautdeg2}
\end{figure}

\section*{Acknowledgements}
Marco Polo Garc\'{\i}a  was partly supported by Conahcyt PhD fellowship grant number 905424.
Ahida Ortiz Santos was partly supported by Conahcyt Master fellowship grant number 1147239.
Martha Alvarez  was partly supported by Programa Especial de Apoyo a Proyectos de Docencia e Investigación  2023, CBI-UAMI.

\bibliographystyle{plain}
\bibliography{ref_plaga}
\end{document}